\documentclass[journal]{IEEEtran}
\IEEEoverridecommandlockouts
\usepackage{comment}
\usepackage{multirow, array}
\usepackage{cite}
\usepackage{eurosym}
\usepackage[mathscr]{euscript}
\usepackage{amsmath,amssymb,amsfonts}
\usepackage{algorithmic}
\usepackage{graphicx}
\usepackage{textcomp}
\usepackage[utf8]{inputenc}
\usepackage[top=2cm, bottom=2cm, left=2cm, right=2cm]{geometry}
\usepackage{float}
\usepackage{caption}
\usepackage{subcaption}
\usepackage{tikz,pgfplots,pgfplotstable,booktabs}
\usepackage{xcolor}
\usepackage{colortbl}
\pgfplotsset{compat=1.15}
\usepackage{comment}
\usepackage{circuitikz}
\usepackage{tikz,pgfplots,pgfplotstable,booktabs}
\usepackage{tikz}
\usepackage{pbox}
\usetikzlibrary{arrows,decorations.pathmorphing,backgrounds,fit,positioning,shapes.symbols,chains, shapes,snakes}
\usepackage{xcolor}
\def\BibTeX{{\rm B\kern-.05em{\sc i\kern-.025em b}\kern-.08em
    T\kern-.1667em\lower.7ex\hbox{E}\kern-.125emX}}
    
\newenvironment{ldescription}[1]
  {\begin{list}{}%
   {\renewcommand\makelabel[1]{##1\hfill}%
   \settowidth\labelwidth{\makelabel{#1}}%
   \setlength\leftmargin{\labelwidth}
   \addtolength\leftmargin{\labelsep}}}
  {\end{list}}
\allowdisplaybreaks
\begin{document}
\title{Cost-driven Screening of Network Constraints for the Unit Commitment Problem}
\author{\'A. Porras, S. Pineda, J. M. Morales, A. Jim\'enez-Cordero.
\thanks{The authors are with the OASYS research group, University of Malaga, Malaga, Spain. E-mail: alvaroporras19@gmail.com; spinedamorente@gmail.com;  juanmi82mg@gmail.com (corresponding author); asuncionjc@uma.es.}
\thanks{This project has received funding in part by the Spanish Ministry of Science and Innovation (AEI/10.13039/501100011033) through project PID2020-115460GB-I00; in part by the European Research Council (ERC) under the European Union's Horizon 2020 research and innovation programme (grant agreement No 755705), and in part by the Junta de Andalucía (JA) and the European Regional Development Fund (FEDER) through the research project P20\_00153. \'A. Porras is also financially supported by the Spanish Ministry of Science, Innovation and Universities through the university teacher training program with fellowship number FPU19/03053. The authors thankfully acknowledge the computer resources, technical expertise and assistance provided by the SCBI (Supercomputing and Bioinformatics) center of the University of Malaga.}}
\maketitle
\begin{abstract}

In an attempt to speed up the solution of the unit commitment (UC) problem, both machine-learning and optimization-based methods have been proposed to lighten the full UC formulation by removing as many superfluous line-flow constraints as possible. While the elimination strategies based on machine learning are fast and typically delete more constraints, they may be over-optimistic and result in infeasible UC solutions. For their part, optimization-based methods seek to identify redundant constraints in the full UC formulation by exploring the feasibility region of an LP-relaxation. In doing so, these methods only get rid of line-flow constraints whose removal leaves the feasibility region of the original UC problem unchanged. In this paper, we propose a procedure to substantially increase the line-flow constraints that are filtered out by optimization-based methods without jeopardizing their appealing ability of preserving feasibility. Our approach is based on tightening the LP-relaxation that the optimization-based method uses with a valid inequality related to the objective function of the UC problem and hence, of an economic nature. The result is that the so strengthened optimization-based method identifies not only redundant line-flow constraints but also inactive ones, thus leading to more reduced UC formulations. 
\end{abstract}
\begin{IEEEkeywords}
	 Unit Commitment, constraint screening, optimization-based method, bounding, cost-driven approach
\end{IEEEkeywords}
\section*{Nomenclature}
The main notation used throughout the text is stated below for  quick  reference.
\subsection{Sets}
\begin{ldescription}{$xxx$}
\item [$\mathcal{N}$] Set of nodes, indexed by $n$.
\item [$\mathcal{L}$] Set of transmission lines, indexed by $l$.
\item [$\mathcal{G}$] Set of generating units, indexed by $g$.
\item [$\mathcal{G}_{n}$] Set of generating units connected to the node $n$.
\item [$\mathcal{T}$] Set of time periods, indexed by $t$.
\end{ldescription}
\subsection{Parameters}
\begin{ldescription}{$xxxxx$}
\item [$a_{ln}$] Power transfer distribution factor (PTDF) of line $l$ at node $n$.
\item [$c_{g}$] Production cost of the generating unit $g$.
\item [\textcolor{black}{$d_{n}$}] {\color{black} Net demand, i.e., load minus renewable generation, at node $n$.}
\item [$\overline{f}_{l}$] Capacity of the transmission line $l$.
\item [$\underline{p}_{g}/\overline{p}_{g}$] Minimum/Maximum output of  unit $g$.
\end{ldescription}
\subsection{Variables}
\begin{ldescription}{$xxx$}
\item [$p_{g}$] Power output dispatch of generating unit $g$.
\item [$q_{n}$] \textcolor{black}{Net injected power at node $n$.}
\item [$u_{g}$] Commitment of generating unit $g$.
\end{ldescription}
\section{Introduction}
\label{sec:introduction}

\subsection{\textcolor{black}{Motivation}}

\IEEEPARstart{T}{he} unit commitment (UC) problem is fundamental in the operation of power systems. It is, actually, at the core of the tools used by Regional Transmission Operators in the US, such as ISO New England and PJM, to clear their day-ahead electricity markets \cite{shahidehpour2003market}. The solution to the UC problem corresponds to the most economically efficient operating schedule, given by the on/off statuses and productions levels of the power plants. The objective of the UC problem is thus to minimize system operation costs, whereas physics and engineering constraints are satisfied \cite{knueven2020mixed}.

Mathematically, the UC problem is typically formulated as a large-scale mixed-integer program (MIP), which belongs to the class of NP-hard problems, even for a single period \cite{bendotti2019complexity}. For this reason, the development of strategies to solve this problem to optimality in a computationally efficient manner has been and is still a popular research topic.

Given the well-known facts that i) dealing with network-constraints in the UC formulation considerably complicates obtaining and certifying the optimal UC plan, and ii) in many systems most transmission lines are oversized and thus are seldom congested, recent research efforts have been put on strategies to get reduced, more compact UC formulations by eliminating superfluous line-flow constraints \cite{yang2021machine}. As pointed out in \cite{pineda2020data}, the reduction process rests on identifying redundant and inactive power flow constraints that can be safely removed from the target UC problem without affecting its optimal solution.

\subsection{\textcolor{black}{Literature review}}

There is a wealth of literature about constraint screening for operational problems in power systems. For instance, spurred by the recent boom of machine learning and artificial intelligence, a number of data-driven strategies have been proposed to detect redundant and inactive constraints very fast by learning from previously solved problem instances. \textcolor{black}{Different strategies to learn the set of active constraints have been applied to the economic dispatch problem \cite{yang2020fast}, the optimal power flow \cite{misra2021learning,deka2019learning,hasan2020hybrid} and the unit commitment problem \cite{pineda2020data,xavier2021learning}.} While computationally inexpensive, the constraint-screening approaches purely based on machine learning carry a risk of constraint misidentification and therefore, may render reduced problems which are not equivalent to the target ones. This means that there exists, by construction, a nonzero probability that the solution to the reduced problem be suboptimal or even infeasible in the original one. \textcolor{black}{These infeasibilities can be removed by combining the machine-learning method with a constraint-generation algorithm, as proposed in [5], at the expense of increasing the computational time.}


On the other hand, the technical literature also includes a family of constraint-screening methods that seek to reduce the UC model as much as possible while ensuring the equivalence between the reduced and original UC problems. For this purpose, these methods usually involve, in one way or another, some form of optimization and consequently, are computationally more demanding in general. Perhaps, the most well-known method within this family is an iterative procedure based on constraint generation, \cite{benameur2006constraint}. In this method, the violated constraints from the original UC problem are gradually added to the reduced one until the solution to the latter is feasible in the former. This procedure has been used, for example, in \cite{xavier2019transmission,tejada2017security} to lighten the computational burden of the SCUC problem by filtering out post-contingency constraints. 
The main disadvantage of the methods based on constraint generation is that they can become computational costly if the number of iterations that are needed to guarantee solution feasibility is too large.

There is another group of optimization-based methods that concentrates on identifying \emph{redundant} constraints, also referred to as \emph{non-umbrella} constraints in \cite{ardakani2013identification}. These are the constraints that,  
if removed, the feasibility region of the original UC problem remains unchanged. The basic idea behind these methods is thus to check whether the constraints that are candidates to be removed are violated or not over a relaxation of this region. If they are not violated, then they can be safely ruled out. In general, this check requires optimization, as it typically translates into solving a series of bounding problems over the relaxed feasible region and, logically, makes sense provided that these bounding problems are much easier to solve than the target optimization problem. \textcolor{black}{Examples of works that follow this \emph{modus operandi} to solve the UC problem are \cite{zhai2010fast, madani2016constraint, roald2019implied, ye2020data}. These bounding-based methods aim to remove as many redundant constraints as possible from the full UC formulation. Nevertheless, even if all the redundant constraints are successfully identified and screened out, there may still remain a number of constraints in the so reduced UC problem that are not needed because they do not oppose to the minimization of the UC cost. In other words, it is the objective function of the UC problem, and not its feasible region, what makes such constraints superfluous. We use the qualifier ``inactive'' to refer to these constraints. Analogous ideas  have been used within the context of the security-constrained DC optimal power flow, see, for instance,   \cite{ardakani2013identification, weinhold2020fast} and references therein.}
\subsection{\textcolor{black}{Contribution}}
In this paper, we propose a strategy to endow bounding-based methods with the ability to detect, and thus discard, not only \emph{redundant} but also \emph{inactive} constraints. \textcolor{black}{On the contrary, the proposed approach retains both \emph{active} and \emph{quasi-active} constraints, as defined in \cite{pineda2020data}.} For this purpose, we build bounding problems that do not only account for technical aspects of the UC problem but also economic ones. To our knowledge, only the recently published paper \cite{ma2021redundant} is in the spirit of our proposal, although they do not follow the trail of bounding-based methods. More specifically, they propose a heuristic algorithm that runs as follows: First, they address the UC without the network, the solution to which they call \emph{cost-based schedule}. Then, they perturb this cost-based schedule by shifting a percentage $\alpha$ of the total generation output (referred to as the \emph{adjustment coefficient})  to the nodes with the highest or lowest \textcolor{black}{Power Transfer Distribution Factors (PTDFs) for a given line}. Finally, they check whether the so perturbed solution congests that line. If it does not become congested, this line is classified as inactive.

{\color{black} Against this background, we develop a bounding approach with the ability to screen out both redundant and inactive constraints. To this end, in the linear programs that compute bounds on the line power flows, we introduce a valid inequality representing a cost budget, i.e., an upper bound on the UC solution cost. In this way, power dispatches that are implausible because of their high cost are made infeasible in the bounding problems. In addition, we combine the cost-budget constraint with a set of plausible nodal net demands that is a convex combination of historically observed values, so that unlikely nodal allocations of the system net demand are discarded. All this together increases the constraint-screening ability of our approach. Finally, we run a series of numerical experiments where we compare the performance of our proposal and its variants with alternative data-driven procedures and bounding techniques for constraint screening.}

The rest of this paper is organized as follows.  Section~\ref{sec:Methodology} elaborates on the proposed procedure and explains the methodology we use to benchmark our approach. Section~\ref{sec:case} discusses simulation results from an illustrative and a more realistic case study. Finally, conclusions are duly drawn in Section~\ref{sec:conclusion}.

\section{Methodology}
\label{sec:Methodology}

\subsection{\textcolor{black}{Unit commitment formulation}}

In this section, we present the proposed methodology and discuss how it compares with existing ones. \textcolor{black}{The vast majority of existing papers that propose methods to screen network constraints consider single-period problems \cite{pineda2020data,yang2020fast,misra2021learning,deka2019learning,hasan2020hybrid,ardakani2013identification,zhai2010fast,madani2016constraint,roald2019implied,weinhold2020fast}. We follow suit and} focus on the single-period UC problem (with no inter-temporal constraints) that is formulated below. 
\begin{subequations}
\label{UC-DCOPF}
\begin{align}
&\min_{u_g,p_g} \hspace{3pt}  \sum_{g \in \mathcal{G}} c_{g} \hspace{1pt} p_{g} \label{eq:UC_objective}\\
&\text{subject to:}\notag\\
& q_n = \sum_{g \in \mathcal{G}_{n}} p_{g} - d_{n}, \quad \forall n \in \mathcal{N} \label{eq:UC_injected}\\
& \sum_{n\in\mathcal{N}} q_n = 0 \label{eq:UC_balance}\\
& u_{g} \hspace{1pt} \underline{p}_{g} \leq p_{g} \leq u_{g} \hspace{1pt} \overline{p}_{g}, \quad \forall g \in \mathcal{G} \label{eq:UC_power}\\
& -\overline{f}_{l} \leq \sum_{n\in\mathcal{N}}a_{ln} q_n \leq \overline{f}_{l}, \quad \forall l \in \mathcal{L} \label{eq:UC_flow}\\
& u_{g} \in \{0, 1\}, \quad \forall g \in \mathcal{G} \label{eq:UC_binary}
\end{align} \label{eq:UC}
\end{subequations}
\noindent The objective function \eqref{eq:UC_objective} minimizes the total production cost. Equation \eqref{eq:UC_injected} computes the net injected power at node $n$, while constraint \eqref{eq:UC_balance} corresponds to the power balance. Constraints \eqref{eq:UC_power} and \eqref{eq:UC_flow} enforce the limits on the power output of the generating unit $g$ and the power flow through the transmission line $l$, respectively. Note that the power flows through the transmission network are modeled using a DC approximation, {\color{black} where $a_{ln}$ in constraint \eqref{eq:UC_flow} stands for the  Power  Transfer  Distribution Factor (PTDF) of line $l$ at node $n$}. Expressions \eqref{eq:UC_binary} set the binary character of the commitment variables.

\subsection{\textcolor{black}{Existing optimization-based screening approaches}} \label{sec:existing}

Problem~\eqref{eq:UC} can be made significantly easier to solve if those constraints~\eqref{eq:UC_flow} with no effect on the optimal unit-commitment plan are removed \cite{ye2020data}. As discussed in the introductory section of this paper, several optimization-based approaches aim at removing redundant line capacity constraints, in particular, \cite{zhai2010fast, roald2019implied}. To do so, these methods check, for each line, whether there is no feasible UC plan under which the line capacity is hit or surpassed. If that is indeed the case, then the corresponding line-flow constraint can be safely screened out because it is essentially redundant. \textcolor{black}{For computational tractability, however, these methods do not work with the feasibility region of the original UC problem, but with that of an LP relaxation  that can be swiftly solved using commercially available optimization software.} In general, for each transmission line $l^{\prime}$, the following two optimization problems are solved:
\begin{subequations}
\label{eq:max_flow}
\begin{align}
&\max_{u_g,p_g,d_n} \hspace{1.5pt} / \hspace{1.5pt} \min_{u_g,p_g,d_n} \hspace{3pt} f_{l^{\prime}} = \sum_{n\in\mathcal{N}}a_{l^{\prime}n} q_n \label{eq:screen_objective}\\
&\text{subject to:}\notag\\
& \eqref{eq:UC_injected}-\textcolor{black}{\eqref{eq:UC_power}}\label{eq:screen_balance}\\
&{\color{black} 0 \leq u_{g} \leq  1,} \quad \forall g \in \mathcal{G} \label{eq:screen_power}\\
& -\overline{f}_{l} \leq \sum_{n\in\mathcal{N}}a_{ln}q_n \leq \overline{f}_{l}, \forall l \in \mathcal{L},l \neq l^{\prime}, \label{eq:screen_flow}\\
%
& \mathbf{d} \in \mathcal{D} \label{eq:screen_demand}
\end{align} 
\end{subequations}

Problem \eqref{eq:max_flow} seeks to maximize/minimize the power flow through line $l^{\prime}$ (denoted as $f_{l^{\prime}}$) over an LP relaxation of the feasible region of the UC problem \eqref{eq:UC}\textcolor{black}{, in which the generated quantities are allowed to take values between zero and the respective generator's capacity $\overline{p}_g$.} Furthermore, the vector of net nodal demands $\mathbf{d}= (d_n)_{n\in \mathcal{N}}$ is turned into a vector of decision variables belonging to the set $\mathcal{D}$. If the minimum and maximum values of $f_{l^{\prime}}$ determined by \eqref{eq:max_flow} lie between  $-\overline{f}_{l^{\prime}}$ and $\overline{f}_{l^{\prime}}$, then the corresponding constraints \eqref{eq:UC_flow} are identified as redundant in problem \eqref{eq:UC} for any~$\mathbf{d} \in \mathcal{D}$. 

Depending on the specific LP relaxation that is used in problem \eqref{eq:max_flow} and the specific set $\mathcal{D}$ that is considered, different methods for identifying redundant line-flow constraints can be derived. For instance, if the line-capacity limits \eqref{eq:screen_flow} are eliminated from \eqref{eq:max_flow} and the set $\mathcal{D}$ is reduced to a singleton, namely, the predicted net demand, then we get the screening method proposed in \cite{zhai2010fast}. Of course, enforcing the power flow constraints~\eqref{eq:screen_flow} in \eqref{eq:max_flow} improves the screening ability of this problem at the expense of increasing its computational burden. Further, if the set $\mathcal{D}$ is such that the net nodal demands $d_n$ are not fixed to a given value but vary within a predefined range, i.e., $\underline{d}_n\leq d_n \leq \overline{d}_n$, then we get the method suggested in \cite{roald2019implied}. The advantage of using a full-dimensional set $\mathcal{D}$ in \eqref{eq:max_flow} is that the line-flow constraints that are filtered out by this problem are, in fact, superfluous for a whole range of operating conditions. 
The downside of this approach is that, in general, the larger the size of $\mathcal{D}$, the lower the number of redundant line-flow constraints that are detected by the collection of problems \eqref{eq:max_flow}. \textcolor{black}{Throughout this paper, we use the method proposed in \cite{roald2019implied} as benchmark (BN).}

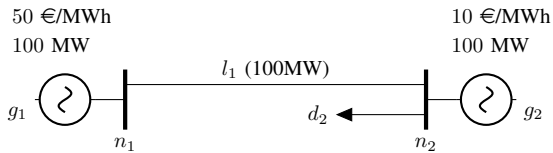
\begin{figure}
\centering
\begin{circuitikz}[scale=0.8,>=triangle 45,every node/.style={scale=0.8}]
\draw [ultra thick] (0,4)  -- (0,3) node[anchor=north]{$n_{1}$};
\draw [ultra thick] (5,4)  -- (5,3) node[anchor=north]{$n_{2}$};
\draw[] (-1.5,3.5) node[anchor=north east]{$g_1$} to [sV] (-0.5,3.5);
\draw (-0.5,3.5) -- (0,3.5);
\draw[] (6.5,3.5) node[anchor=north west]{$g_2$} to [sV] (5.5,3.5); 
\draw (5.5,3.5) -- (5,3.5); 
\draw (0,3.75) -- (5,3.75);
\node[] at (2.5,4) {$l_1$ (100MW)};
\draw[->] (5,3.25) -- (3.5,3.25) ;
\node[anchor=east] at (3.5,3.25) {$d_2$};
\node[anchor=west] at (-2,4.9) {$50$ \euro/MWh};
\node[anchor=west] at (-2,4.4) {$100$ MW};
\node[anchor=west] at (5.3,4.9) {$10$ \euro/MWh};
\node[anchor=west] at (5.3,4.4) {$100$ MW};
\end{circuitikz}
\caption{Two-node illustrative example}
\label{fig:2bus}
\end{figure}

A relevant drawback of LP relaxation \eqref{eq:max_flow} is that it only accounts for aspects of the UC problem that are purely technical. That is, it completely ignores the economics behind the UC problem, which is embodied in the minimization of the production costs. As a result, the series of problems~\eqref{eq:max_flow} are only able to screen out \emph{redundant} constraints \textcolor{black}{but fail to identify constraints that are \emph{inactive} at the optimum.}
To illustrate this, consider the two-node network depicted in Fig. \ref{fig:2bus} with an expensive generator in $n_1$, a cheap generator in $n_2$, a line with a capacity of 100 MW, and a net load $d_2$ that varies between 80 MW to 120 MW. Note that line $l_1$ never becomes congested since $d_2$ would be first satisfied by the cheaper unit connected to the same node and the maximum production of the expensive unit would be 20 MW without any network congestion. However, if problem \eqref{eq:max_flow} is solved for $l_1$ with $\mathcal{D}$ being the interval $80\leq d_2 \leq 120$, the possible maximum power flow through this line would reach 100 MW by fully dispatching $g_1$ and therefore, the capacity limit constraint would not be removed. In fact, such a constraint would be kept regardless of the marginal generating cost assigned to each of the two units. 

\subsection{\textcolor{black}{Proposed cost-driven constraint screening}} \label{sec:proposed}

To overcome the drawback of existing optimization-based screening approaches and reduce even further the UC problem~\eqref{eq:UC}, we propose to tighten the LP relaxation used in the set of problems \eqref{eq:max_flow} by including a valid inequality on the optimal objective function of~\eqref{eq:UC}. More specifically, we propose to find the maximum power flow through each line $l'$ by solving \eqref{eq:max_flow_ub}, where constraint \eqref{eq:maximum_cost} imposes an upper bound on the total production cost. Thus, inactive but non-redundant constraints can  also be screened out depending on the cost of the unit-commitment plan. This methodology is denoted as UB in this paper.
\begin{subequations}
\label{eq:max_flow_ub}
\begin{align}
&\max_{u_g,p_g,d_n,q_n} \hspace{1.5pt} / \hspace{1.5pt} \min_{u_g,p_g,d_n,q_n} \hspace{3pt} f_{l^{\prime}} = \sum_{n\in\mathcal{N}}a_{l^{\prime}n}q_n \\
&\text{subject to:}\notag\\
& \eqref{eq:screen_balance} - \eqref{eq:screen_demand}\\
& \sum_{g \in \mathcal{G}} c_{g} \hspace{1pt} p_{g} \leq \overline{C} \label{eq:maximum_cost}
\end{align} 
\end{subequations}
Coming back to the illustrative example in Section \ref{sec:existing}, if $d_2$ never exceeds 120 MW, then the maximum generation cost would be equal to \euro2000. Solving problem \eqref{eq:max_flow_ub} for $l_1$, the range $80\leq d_2 \leq 120$, and constraint \eqref{eq:maximum_cost} formulated as $50p_{g_1}+10p_{g_2}\leq2000$ yields a maximum power flow through that line below its capacity, since $g_1$ cannot be dispatched above 40 MW. Consequently, the capacity constraint of this line can be safely removed for these particular values of the units' marginal production costs.

\begin{figure}
    \centering
    \begin{tikzpicture}[scale=0.8]
	\begin{axis}[xmin=65,xmax=70,ymin=755,ymax=930,
	xlabel=Aggregate Net Demand,
	ylabel=Operating Cost,
	ytick=\empty,
	xtick={65.09,68.164,69.710}, 
	xticklabels={$\underline{D}_1$,$\overline{D}_1=\underline{D}_2$,$\overline{D}_2$},ticklabel style = {}]
	\addplot[only marks, gray, mark=*,mark options={scale=0.4, fill=gray},text mark as node=true] table[x=AD,y=Cost,col sep=comma] {QSLR.csv};
	\addplot[thick] coordinates {(65.09,762.195)(68.164,850.926)};
	\addplot[thick] coordinates {(68.164,856.6)(69.710,922.782)};
	\addplot[thick,dashed] coordinates {(65,911.297)(69.583,911.297)};
	\addplot[thick,dotted] coordinates {(65.09,755)(65.09,762.195)};
	\addplot[thick,dotted] coordinates {(68.164,755)(68.164,856.6)};
	\addplot[thick,dotted] coordinates {(69.710,755)(69.710,922.782)};
    \end{axis}
    \node[anchor=south,sloped] at (-0.5,4.9) {$C^{\max}$};
\end{tikzpicture}
\caption{$1$-quantile piecewise linear regression}
\label{fig:ilus_approach}
\end{figure}
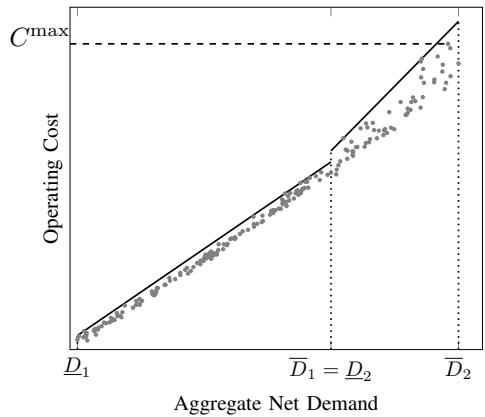

Taken from a more realistic system, Fig. \ref{fig:ilus_approach} displays an example of the total operating cost as a function of the aggregated net demand for 200 past instances of problem \eqref{eq:UC}, with $C^{\max}$ being the highest observed cost. As expected, similar aggregated net demand values may lead to different operating costs due to distinct allocations of such a demand among network buses. Obviously, fixing $\overline{C}=C^{\max}$ in \eqref{eq:max_flow_ub} would be a too conservative strategy. Indeed, such an upper bound would be much higher than the actual operating cost for low net demand values and therefore, constraint~\eqref{eq:maximum_cost} would be too loose and useless for our purpose. In other words, the constraint-screening ability of problems~\eqref{eq:max_flow_ub} would not be any better than that of \eqref{eq:max_flow} in this case. Accordingly, to take full advantage of the valid inequality \eqref{eq:maximum_cost}, we make  $\overline{C}$ dependent on the net demand vector $\mathbf{d}$, i.e. $\overline{C} = \overline{C}(\mathbf{d})$. In doing so, tighter upper bounds can be obtained for all net demand levels and a larger number of line capacity constraints can be screened out.

We justify next our modeling choices to estimate the function $\overline{C}(\mathbf{d})$. First, to reduce model complexity and the risk of overfitting, we make the upper-bound $\overline{C}$ dependent on the aggregated net demand $D=\sum_{n\in\mathcal{N}}{d_n}$ only, as opposed to using all nodal net demands as explanatory variables. \textcolor{black}{Nevertheless, for systems with a high penetration of renewable generation, the operating cost may be better explained by using two explanatory variables, namely, the aggregate demand and the aggregate renewable generation.} Second, in order to keep problems \eqref{eq:max_flow_ub} manageable and computationally tractable, we approximate the relation between the upper-bound $\overline{C}$ and the aggregate net demand $D$ through a piecewise linear function. Third, we resort to  $\tau$-quantile regression, with $0\leq\tau\leq 1$,  \cite{koenker2001quantile}. More particularly, for each segment, we use $1$-quantile regression in order to ensure that the approximated upper bound is higher than the operating cost of all observed instances. The solid line in Fig.~\ref{fig:ilus_approach} represents the $1$-quantile piecewise linear function for this illustrative data set. 

The quantile piecewise linear regression we propose is then characterized by a set of segments $\mathcal{S}$, indexed by $s$. \textcolor{black}{We can resort to well-established procedures in Statistics to decide on the number of segments. For instance, we can borrow the popular \emph{elbow criterion} from the realm of data clustering to this end. According to this criterion, we compute and plot the quantile loss as a function of the number of segments and the \emph{elbow} of the curve is taken as a good value for the number of pieces to be picked. The idea is to choose the number of segments from which the implied reduction in the quantile loss starts to diminish rapidly. The choice of the number of segments could also be based on even simpler procedures such as the visual inspection of the data scatter plot ``net demand vs. cost'' (see Fig. 2) and/or more sophisticated tuning techniques based on bootstrapping and cross-validation.}

The minimum and maximum aggregate net demand for each segment is denoted by $\underline{D}_s$ and $\overline{D}_s$, and we introduce a binary variable $y_s$ per segment, which is equal to 1 if $\underline{D}_s \leq D \leq \overline{D}_s$, and to 0 otherwise. For illustration, Fig. \ref{fig:ilus_approach} includes $\underline{D}_s$ and $\overline{D}_s$ for the two segments considered. The upper bound for an aggregate net demand $D$ can be thus computed as:
\begin{equation*}
    \overline{C}(D) = \sum_{s \in \mathcal{S}} y_s(\rho_s + \nu_s D)
\end{equation*}
where $\rho_s$ and $\nu_s$ are the intercept and slope of the linear function corresponding to segment $s$. Under this approximation, problem \eqref{eq:max_flow_ub} can be recast as:
\begin{subequations}
\label{eq:max_flow_regression}
\begin{align}
&\max_{u_g,p_g,d_n,q_n,y_s,D} \hspace{1.5pt} / \hspace{1.5pt} \min_{u_g,p_g,d_n,q_n,y_s,D} \hspace{3pt} f_{l^{\prime}} = \sum_{n\in\mathcal{N}}a_{l^{\prime}n}q_n \\
&\text{subject to:}\notag\\
& \eqref{eq:screen_balance} - \eqref{eq:screen_demand}\\
& \sum_{g \in \mathcal{G}} c_{g} p_{g} \leq \sum_{s \in \mathcal{S}} y_{s}\left( \rho_{s} + \nu_{s} \, D \right)  \label{eq:upperbound_quantile}\\
& D = \sum_{n \in \mathcal{N}} d_{n}\label{eq:aggreagte_demand}\\
& \sum_{s \in \mathcal{S}} y_{s}\underline{D}_{s} \leq D \leq \sum_{s \in \mathcal{S}} y_{s}\overline{D}_{s}\label{eq:aggreagte_demand_bounds}\\
& \sum_{s \in \mathcal{S}} y_{s} = 1 \label{eq:activate_segment}\\
& y_{s} \in \{0, 1\}, \quad \forall s \in \mathcal{S}. \label{eq:binary_character_y}
\end{align} 
\end{subequations}  
Constraint \eqref{eq:upperbound_quantile} is identical to \eqref{eq:maximum_cost}, except that the constant upper bound is replaced with the piecewise linear function obtained from the $1$-quantile regression. Equation \eqref{eq:aggreagte_demand} computes the aggregate net demand $D$. Constraints \eqref{eq:aggreagte_demand_bounds} impose segment-dependent minimum and maximum bounds on the aggregate net demand, in that order. Therefore, if $y_s=1$, then $\underline{D}_s\leq D \leq \overline{D}_s$. Expression \eqref{eq:activate_segment} ensures a one-to-one mapping between the aggregate net demand segment and the cost piece. Finally, constraints \eqref{eq:binary_character_y} set the binary character of the decision variables $y_s$. \textcolor{black}{The product of binary and continuous variables, $y_s D$ in \eqref{eq:upperbound_quantile}, can be easily linearized by using integer algebra results \cite{floudas1995} so that \eqref{eq:max_flow_regression} can be formulated as a mixed-integer linear optimization problem and solved using commercial optimization software. }


A second drawback of LP relaxation \eqref{eq:max_flow} relates to the choice of set $\mathcal{D}$. For instance, the method proposed in \cite{roald2019implied} defines $\mathcal{D}$ as the Cartesian product of the intervals $[\underline{d}_n, \overline{d}_n]$. Under this assumption, the solution of problems \eqref{eq:max_flow} may render an implausible net demand profile $\mathbf{d}= (d_n)_{n\in \mathcal{N}}$, which does not conform at all with the observed spatial correlations among the net nodal loads. For this reason, in this paper, we also explore defining $\mathcal{D}$ as the set of net demand profiles $\mathbf{d}= (d_n)_{n\in \mathcal{N}}$ that are a convex combination of the observed ones. This is imposed as follows:
\begin{subequations}
\label{eq:convex_combination}
\begin{align}
& d_n = \sum_{t\in\mathcal{T}} \alpha_{t} \tilde{d}_{tn} \label{eq:convex_combination_1}\\
& \alpha_{t} \geq 0, \quad \textcolor{black}{\forall t\in \mathcal{T}} \label{eq:convex_combination_2}\\
& \sum_{t\in\mathcal{T}} \alpha_{t} = 1, \label{eq:convex_combination_3}
\end{align}
\end{subequations}
\noindent where $\tilde{d}_{tn}$ denotes the net demand at node $n$ in past time period $t$. {\color{black}We remark that we may allow the sum of the scalars $\alpha_t$ in (5c) to be equal  to a number greater than one, if we want to increase the probability of not making a mistake when removing a particular line-flow constraint. However, based on the rationale that demand patterns do not change quickly, our choice of $\mathcal{D}$ as (5) offers a good compromise  between risk of mistake and constraint-screening power of the bounding problem (2)  for UC instances that are to be solved in the near future. Furthermore, the convex hull of past net demand observations we use as $\mathcal{D}$ can be frequently enriched with new data points as we move forward in time.}

\textcolor{black}{Naturally, the performance of the proposed approach requires that the set $\mathcal{T}$ includes a sufficiently high number of past time periods that characterize the variability of the system operating conditions.} Defining $\mathcal{D}$ through 
\eqref{eq:convex_combination} instead of as the Cartesian product of the intervals $[\underline{d}_n, \overline{d}_n]$ results in a more constrained problem \eqref{eq:max_flow} and a higher number of removed constraints. \textcolor{black}{In this paper, this modeling choice is denoted as CC (which stands for \emph{Convex Combination}).} 

In summary, methods UB and CC shrink the feasible space of generation and net demand variables in problems \eqref{eq:max_flow}, respectively, to leave out unrealistic operating conditions. We also advocate for the synergistic use of UB and CC, that is, for the constraint screening method UB+CC, which exploits an upper-bound function on the total operating cost over a convex combination of the net demand profiles. For clarity, Table \ref{tab:approaches} compiles the  optimization problems solved by each of the approaches that are considered in this paper.


\subsection{Performance evaluation procedure}\label{sec:Eval}

The procedure to measure the performance of the approaches described in Section \ref{sec:proposed} to screen out line capacity constraints and reduce the computational burden of the UC problem runs as follows:
\begin{itemize}
    \item[1)] Use historical information on past unit commitment instances to adjust some of the parameters of the screening optimization problems. For instance, in BN, this information can be used to adjust the range of each nodal net demand that defines the set $\mathcal{D}$, while, in UB, this data is employed to estimate the parameters of the 1-quantile piecewise linear regression. Similarly, imposing the convex combination of net demands in CC requires past net demand profiles. Finally, UB+CC uses historical information to compute both the regression parameters and convex combinations of demands.
    \item[2)] Depending on the approach, solve the optimization problem listed in Table \ref{tab:approaches} for each of the lines of the network and determine the set of line capacity constraints that can be removed from the original UC problem. In particular, if the maximum (minimum) flow of line $l$ determined by the screening optimization problem is below (above) the limit $\overline{f}_l$ ($-\overline{f}_l$), then constraint $f_l \leq \overline{f}_l$ ($-\overline{f}_l \leq f_l $) can be screened out.  
    \item[3)] Solve a reduced unit commitment problem similar to \eqref{eq:UC} without including the line capacity constraints \eqref{eq:UC_flow} that have been filtered out in step 2).
    \item[4)] Fix the binary commitment decisions to those obtained in step 3) and solve the unit commitment problem including all constraints. 
    \item[5)] \textcolor{black}{Measure the performance of the screening approach in terms of i) the average number of line capacity constraints that are filtered out in step 2), ii) the average computational time required to solve the reduced UC problem of step 3) with respect to the full UC formulation, and iii) the average optimality loss of the obtained solution with respect to that of the full UC problem.}
\end{itemize}

\begin{table}[]
\caption{Constraint screening approaches}
\begin{center}
\begin{tabular}{ll}
\hline
Approach & Screening optimization problem \\
\hline
BN & \eqref{eq:screen_objective} s.t. \eqref{eq:screen_balance} - \eqref{eq:screen_demand}   \\
UB & \eqref{eq:screen_objective} s.t. \eqref{eq:screen_balance} - \eqref{eq:screen_demand}, \eqref{eq:upperbound_quantile}-\eqref{eq:binary_character_y}   \\
CC & \eqref{eq:screen_objective} s.t. \eqref{eq:screen_balance} - \eqref{eq:screen_demand}, \eqref{eq:convex_combination_1}-\eqref{eq:convex_combination_3}   \\
UB+CC & \eqref{eq:screen_objective} s.t. \eqref{eq:screen_balance} - \eqref{eq:screen_demand}, \eqref{eq:upperbound_quantile}-\eqref{eq:binary_character_y}, \eqref{eq:convex_combination_1}-\eqref{eq:convex_combination_3} \\
\hline 
\end{tabular}  
\label{tab:approaches}
\end{center}
\end{table}

\section{Numerical Results}\label{sec:case}

This section presents numerical results on two power systems. The first one is a small-size system especially designed for illustrative purposes. The second one is based on a realistic power system from Texas. 
All experiments have been carried out on a cluster with 21 Tb RAM, running Suse Leap 42 Linux distribution. The constraint screening approaches presented in Table \ref{tab:approaches} have been modeled using Pyomo 5.7.1 \cite{pyomo}  and solved with CPLEX 20.1 \cite{cplex}. \textcolor{black}{The number of segments of the $1$-quantile piecewise linear regression has been set to one in Section \ref{subs: Illustrative Example} (illustrative example) and to three in Section \ref{subs: Results for the Large-scale Case} (realistic case study). In the latter case, our choice is motivated by the visual inspection of the data scatter plot \emph{net demand vs. cost}.} The breakpoints have been found with the Python library \texttt{PWLF}.



\subsection{Illustrative Example}\label{subs: Illustrative Example}

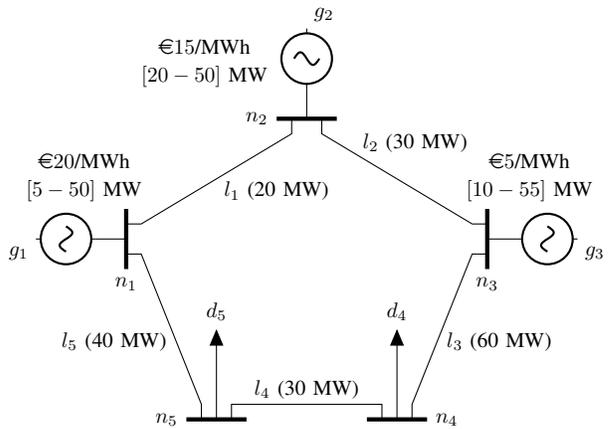
\begin{figure}
\centering
\begin{circuitikz}[scale=0.8,>=triangle 45,every node/.style={scale=0.8}] 
\draw [ultra thick] (0,2.5)  -- (0,3.5) node[anchor=north] at (0,2.5) {$n_{1}$};
\draw [ultra thick] (2.5,5)  -- (3.5,5) node[anchor=east] at (2.5,5) {$n_{2}$};
\draw [ultra thick] (6,2.5) -- (6,3.5) node[anchor=north] at (6,2.5) {$n_{3}$};
\draw [ultra thick] (4,0)  -- (5,0) node[anchor=west] at (5,0) {$n_{4}$};
\draw [ultra thick] (1,0)  -- (2,0) node[anchor=east] at (1,0) {$n_{5}$};
\draw (0,3.25) -- (0.25,3.25) -- (2.75,4.75) -- (2.75,5) node[] at (2.5,3.8) {$l_{1}$ (20 MW)};
\draw (3.25,5) -- (3.25,4.75) -- (5.75,3.25) -- (6,3.25) node[] at (4.8,4.6) {$l_{2}$ (30 MW)};
\draw (6,2.75) -- (5.75,2.75) -- (4.75,0.25) -- (4.75,0)
node[] at (6.2,1.3) {$l_{3}$ (60 MW)};
\draw (4.25,0) -- (4.25,0.25) -- (1.75,0.25) -- (1.75,0) node[] at (3,0.5) {$l_{4}$ (30 MW)};
\draw (1.25,0) -- (1.25,0.25) -- (0.25,2.75) -- (0,2.75) node[] at (-0.2,1.3) {$l_{5}$ (40 MW)};
\draw[] (-1.5,3) node[anchor=north east]{$g_1$} to [sV] (-0.5,3); 
\draw (-0.5,3) -- (0,3); 
\draw[] (3,6.5) node[anchor=south west] {$g_2$} to [sV] (3,5.5);
\draw (3,5.5) -- (3,5);
\draw[] (7.5,3) node[anchor=north west]{$g_3$} to [sV] (6.5,3); 
\draw (6.5,3) -- (6,3); 
\draw[->] (1.5,0) -- +(0,1.5); 
\node[] at (1.5,1.8) {$d_5$};
\draw[->] (4.5,0) -- +(0,1.5); 
\node[] at (4.5,1.8) {$d_4$};
\node[] at (-0.7,4.3) {\euro$20$/MWh};
\node[] at (-0.7,3.8) {$[5-50]$ MW};
\node[] at (1.3,6.2) {\euro$15$/MWh};
\node[] at (1.3,5.7) {$[20-50]$ MW};
\node[] at (6.7,4.3) {\euro$5$/MWh};
\node[] at (6.7,3.8) {$[10-55]$ MW};
\end{circuitikz}
\caption{Five-node system}
\label{fig:5bus}
\end{figure}

Next we illustrate the most salient features of  our approach using the five-node system represented in Figure~\ref{fig:5bus}. This small system includes three thermal units whose linear operating costs and minimum/maximum power limits are indicated in the figure. For simplicity, the five lines in this small system have the same susceptance of 1 p.u., while their capacity limits are also provided in said figure. Finally, Table~\ref{tab:Data-IE} contains historical unit commitment information including the net demand profile, the optimal operating cost and the congested lines for three times periods. 

Due to the capacity limit of line $l_4$, which connects nodes $n_4$ and $n_5$, operating costs may significantly change depending on how the total net demand is distributed among $d_4$ and $d_5$. For instance, in periods $t_1$ and $t_2$, all net demand is located at $n_4$ and can thus be satisfied by the cheaper units $g_2$ and $g_3$ without involving network congestion. Conversely, in period $t_3$ all net demand is at $n_5$, which leads to the congestion of $l_4$ and a higher operating cost due to the commitment of the expensive unit $g_1$. 

\begin{table}
\caption{Historical data -- Illustrative Example}
\vspace{-0.2cm}
\begin{center}
\begin{tabular}{ccccc}
\hline
Time period &$d_{4}$ (MW) &$d_{5}$ (MW)  &  Cost (\euro) & Congested lines\\
\hline
$t_1$  &55    &0   &275.0    &- \\
$t_2$  &75    &0   &575.0    &- \\
$t_3$  &0   &69  &772.5  &$l_4$ \\
\hline
\end{tabular}
\label{tab:Data-IE}
\end{center}
\end{table}

Now consider a new time period $t_4$ with a demand profile of $(d_4,d_5)=(58, 13.8)$\,MW. Solving the complete unit commitment \eqref{eq:UC} for period $t_4$ yields an operating cost of \euro611. For such a demand profile, the power flow through $l_1$ is 7 MW, which is below its limit and thus, the capacity constraint of this line is not binding. However, if such a constraint were removed, the obtained solution would become infeasible since the power flow through $l_1$ would end up exceeding 20 MW. {\color{black} This type of constraint is referred to as \emph{quasi-active} in \cite{pineda2020data}.}
Therefore, to obtain the same solution as the full unit commitment problem, the line capacity constraint of $l_1$ must be enforced.  

%



\begin{table}
\caption{Retained capacity constraints in the UC problem at time period $t_4$ -- Illustrative Example}
\vspace{-0.2cm}
\begin{center}
\renewcommand{\arraystretch}{1.5}
\begin{tabular}{lcccc}
\hline
 & BN & UB & CC & UB+CC \\ 
\hline
$-\overline{f}_{l_1} \leq f_{l_1}$ & x &  &  & \\
$\phantom{-}f_{l_1} \leq \overline{f}_{l_1} $ & x & x & x & x \\
$-\overline{f}_{l_2} \leq f_{l_2}$ &  &  &  &  \\
$\phantom{-}f_{l_2} \leq \overline{f}_{l_2} $ & x & x & x & \\
$-\overline{f}_{l_3} \leq f_{l_3}$ &  &  &  &  \\
$\phantom{-}f_{l_3} \leq \overline{f}_{l_3} $ & x & x &  &  \\
$-\overline{f}_{l_4} \leq f_{l_4}$ & x & x & x & x \\
$\phantom{-}f_{l_4} \leq \overline{f}_{l_4} $ & x & x & x & x \\
$-\overline{f}_{l_5} \leq f_{l_5}$ & x & x & x &  \\
$\phantom{-}f_{l_5} \leq \overline{f}_{l_5} $ &  &  &  &  \\
\hline
\# Retained constraints & 7 & 6 & 5 & 3 \\
\hline
\end{tabular}
\label{tab:Results-IE}
\end{center}
\end{table}

Table \ref{tab:Results-IE} shows the line capacity constraints that are retained in the unit commitment problem for time period $t_4$ using the different approaches described in Section \ref{sec:Methodology}.
For the BN approach, problems \eqref{eq:max_flow} are solved for the five lines considering that demands can vary between the minimum and maximum observed values, that is, 
\begin{subequations}
\label{eq:setD_example}
\begin{align}
  & 0\leq d_4\leq 75 \\
  & 0\leq d_5\leq 69
\end{align}
\end{subequations}
The results obtained by the BN approach indicate that three out of the ten line capacity constraints are redundant and therefore, the reduced unit commitment problem includes seven line capacity constraints.

\begin{figure}
    \centering
    \begin{tikzpicture}[scale=0.8]
	\begin{axis}[xmin=50,xmax=80,ymin=50,ymax=1200,xlabel=Aggregate Net Demand (MW),ylabel=Operating Cost (\euro{}),ticklabel style = {font=\tiny}]
	\addplot[thick] coordinates {(50,97.32)(80,1163.39)};
	\addplot[thick,dashed] coordinates {(71.8,50)(71.8,872)};
	\node[] at (axis cs: 55,200) {$t_1$};
	\node[circle,fill=black,scale=0.6] at (axis cs: 55,275) {};
	\node[] at (axis cs: 69,700) {$t_3$};
	\node[circle,fill=black,scale=0.6] at (axis cs: 69,772.5) {};
	\node[] at (axis cs: 73,560) {$t_4$};
	\node[circle,fill=gray,scale=0.6] at (axis cs: 71.8,611) {};
	\node[] at (axis cs: 76,510) {$t_2$};
	\node[circle,fill=black,scale=0.6] at (axis cs: 75,575) {};
    \end{axis}
\end{tikzpicture}
    \caption{Upper-bound Estimation -- Illustrative Example}
    \label{fig:estimation}
\end{figure}
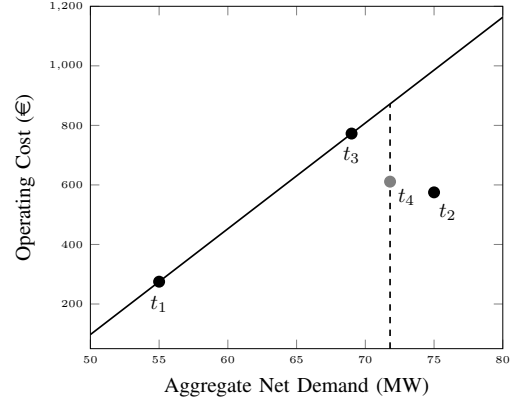

To illustrate the UB approach, Fig. \ref{fig:estimation} plots the operating cost as a function of the aggregate net demand of time periods $t_1$, $t_2$ and $t_3$. The continuous line represents the 1-quantile regression described in Section \ref{sec:Methodology} for a single segment, whose expression is $\overline{C}(D) = -1679.5+35.5D$.
This linear function allows us to estimate an upper bound on the operating cost for new values of the aggregate net demand as follows:
\begin{equation}
    20p_{g_1} + 15p_{g_2} + 5p_{g_3} \leq -1679.5 + 35.5(d_4+d_5) \label{eq:ub_example}
\end{equation}
The upper bound corresponding to $t_4$ is equal to \euro872, which is higher than the actual cost of \euro611. Considering such an upper bound on the operating cost avoids that the power flow through line $l_1$ reaches its minimum value and therefore,  UB yields a reduced problem that includes six line capacity constraints.

Following the CC approach, the set $\mathcal{D}$ formulated in \eqref{eq:setD_example} is replaced by \eqref{eq:setD_example_convex} so that $\mathbf{d}= (d_4, d_5)$ is guaranteed to be a convex combination of the data in Table \ref{tab:Data-IE}. In this case, CC retains the five line capacity constraints indicated in Table~\ref{tab:Results-IE}. 
%
\begin{subequations}
\label{eq:setD_example_convex}
\begin{align}
    & d_4 = \alpha_1 \cdot 55 + \alpha_2 \cdot 75 + \alpha_3 \cdot 0 \\
    & d_5 = \alpha_1 \cdot 0 + \alpha_2 \cdot 0 + \alpha_3 \cdot 69 \\
    & \alpha_1 + \alpha_2 + \alpha_3 = 1
\end{align} \label{eq:convex_example}
\end{subequations}
Finally, the UB+CC approach screens out line capacity constraints by considering both the upper bound on the operating cost computed in \eqref{eq:ub_example} and the convex combination of net demand profiles formulated in \eqref{eq:convex_example}. In doing so, this approach is able to successfully restrict the search space of problems \eqref{eq:max_flow} so that it takes into account information on both the feasible region and the objective function of the UC problem. This way, the UB+CC approach states that only three line capacity limits are to be retained in the unit commitment problem. 
%

The four approaches compared in this illustrative example correctly retain the maximum capacity limit of line $l_1$ and hence, all of them obtain the same UC solution as that of the full UC problem. On the other hand, the number of retained constraints significantly vary between the different approaches. The BN method proposed in the literature is the most conservative one and only removes three constraints. The other three approaches screen out a higher number of constraints by solving problems \eqref{eq:max_flow} with a feasible region that has been tightened using a cost budget and/or a shrunk net demand set $\mathcal{D}$. The UB+CC approach  screens out the largest number of constraints (seven out of ten), which would result in the highest computational savings.



\subsection{Realistic Case Study} \label{subs: Results for the Large-scale Case}
In this section, we consider a realistic power system from Texas, \cite{2000bus}. The network consists of $2000$ buses and $3206$ lines. Hence, the associated UC problem includes $6412$ line capacity constraints. To work with more challenging UC instances, the minimum generating capacity limit, $\underline{p}_g$, has been set as the maximum between the nominal value given in \cite{2000bus} and the $10\%$ of the unit's capacity, as done in \cite{pineda2020data,roald2019implied}. In addition, the marginal costs of the generators $c_g$ have been modified so that they vary in the range \euro$[0, 40]$/MWh. The line capacities, $\overline{f}_l$, are also reduced to have a more congested network. Finally, we synthetically generate the nodal net demands for 8640 time periods as follows. First, a demand value $L$ for each period is randomly sampled from an uniform distribution in the range $[50, 70]$ GW. Second, the net demand at node $n$ is also randomly sampled from a uniform distribution in the range $[0.95L\xi_n, 1.05L\xi_n]$, where $\xi_n$ are the nodal allocation factors of the total system net load taken from \cite{2000bus}. All the data used in this case study is available for download at \cite{Screenpowersystems2021}. {\color{black} Of the 8640 values of nodal net loads we produce, 7200 are used for training and parameter fitting, and the remaining 1440 for testing and performance comparison. In particular, the training data is employed to conduct the quantile regression on which method UB relies.

Our numerical evaluation consists of {\color{black} four} different, but related experiments, which we describe and discuss in the following.

\subsubsection{Experiment 1 -- Fixed load}

Firstly, we consider the setup in which the demand set $\mathcal{D}$ in \eqref{eq:screen_demand} is reduced to a singleton. Specifically, the only element in $\mathcal{D}$ is the vector of nodal net demands in the hour in question (which is assumed to be known or predicted). In this way, we can compare the bounding-based constraint-screening methods BN and UB with the data-driven and heuristic ones proposed in \cite{pineda2020data} and \cite{ma2021redundant}, respectively. Notice that, when $\mathcal{D}$ is a singleton, methods CC and UB+CC become correspondingly equivalent to BN and UB. In the case of the data-driven constraint-screening method in \cite{pineda2020data}, here we report results for 50 neighbors only (which is why we have denoted it as DD50). Nevertheless, similar empirical findings are obtained for a number of neighbors equal to 5 and 500. As for the heuristic constraint-screening method introduced in~\cite{ma2021redundant}, which we denote as H, the tuning parameter on which this method depends has been set to 3\% as the authors in~\cite{ma2021redundant} suggest. 

The main results from this first experiment are collated in Table~\ref{tab:Fixed}. The heuristic method turns out to be the most conservative one, in the sense that it is the one retaining the largest number of line-flow constraints in the reduced UC formulation. Consequently, it leads to the highest solution time. At the other end of the spectrum is the data-driven constraint-screening method DD50, which only keeps 4.8\% of the line-flow limits, thus featuring the lowest computational burden by far. However, its remarkable speed comes at a significant cost. Indeed, the daring nature of DD50 results in UC plans that are infeasible in a few instances and suboptimal in many of them, with an average cost error of 0.022\%. From among the four methods analyzed in Table~\ref{tab:Fixed}, UB is the one showing a more balanced performance. In fact, the cost-based valid inequality this method utilizes gets to reduce the number of retained constraints in more than eight percentage points as compared to BN, which, in turn, means a two-times speedup with respect to this one. Most importantly, these improvements come without delivering infeasible or suboptimal UC plans.

\begin{table}[]
\caption{Results -- Fixed load}
\centering
\begin{tabular}{lcccccc}
\hline
& DD50 & BN & UB & H\\
\hline
Retained constraints (\%) &4.8 &18.5 &10.2 &23.5 \\
\#Infeasibilities &4 &0 &0 &0\\
\#Sub-optimal solutions &540 &0 &0 &0 \\
Cost error (\%) &0.022 &0.000 &0.000 &0.000\\
Screening time (s) &0.05 &0.15 &0.18 &0.02 \\
Reduced UC time (s) &0.39 &3.10 &1.50 &3.34 \\
Total time (s) &0.45 &3.25 &1.68 &3.36 \\
Computational burden (\%) &3.9 &28.2 &14.5 &29.1\\
\hline
\end{tabular}
\label{tab:Fixed}
\end{table}

 } 
 
\subsubsection{Experiment 2 -- Base Case}
{\color{black}
We now consider the alternative setup in which the net demand set $\mathcal{D}$ is not a singleton, but spans a subset of $\mathbb{R}^{|\mathcal{N}|}$. In the case of BN and UB, this subset is given by the Cartesian product of the intervals $[\underline{d}_n, \overline{d}_n]$, $n \in \mathcal{N}$, where the bounds $\underline{d}_n$ and $\overline{d}_n$ are respectively taken as the minimum and maximum values of net load at node $n$ that are observed in the training data set. In the case of CC and UB+CC, the set $\mathcal{D}$ corresponds to the convex combination of the net demand vectors that make up the same training data. This experiment excludes the constraint-screening methods DD50 and H, because these are designed to work with a fixed net load only. In effect, these methods are not intended to filter out line-flow constraints that are superfluous \emph{for a (wide) range of net load values}. Consequently, unlike BN, UB, CC and UB+CC, methods DD50 and H are to be rerun every hour or every day, that is, every time the predicted vector of nodal net loads is updated. 

}

The first important observation to be made is that the three proposed methods (UB, CC, UB+CC) provide, for the 1440 time periods in the test set, the same optimal UC cost as that obtained by the benchmark approach. Consequently, the performance of these methods, {\color{black} which is summarized in Table \ref{tab:Large}, is assessed and compared} in terms of the average percentage of retained constraints and the average computational burden relative to the computational time required to solve the full UC problem. \textcolor{black}{This computational burden does not include the time of the screening procedure, which is assumed to be run offline only once, since the network constraints that are screened out by these approaches can be left out of the UC problem for a wide range of operating conditions. The average computational times of the screening problems to be solved under methods BN, UB, CC and UB+CC amount to 2.9s, 9.7s, 122.4s, and 350.4s, respectively.}




Results in Table \ref{tab:Large} show that method UB involves slight reductions in both the number of retained constraints and the computational burden. Imposing a convex combination of the nodal demand through method CC reduces the retained constraints and computational time by 9\% and 29\%, respectively. Relevantly, it is the synergistic effect of UB and CC, i.e., the proposed approach UB+CC, which involves the largest improvement. Indeed, this method obtains reduced UC problems with 15\% fewer constraints than BN with a computational time 45\% lower on average. As pointed out in Section~\ref{sec:Methodology}, these results demonstrate that the estimated cost-budget function $\overline{C}(D)$ is most effective when combined with a consistent feasible demand set $\mathcal{D}$, like the one used in CC. 

{\color{black} Despite the fact that there are no infeasibilities or suboptimal solutions, these could appear when the cost budget function is underestimated or the unseen net nodal demand vector $\mathbf{d}$ does not belong to the set $\mathcal{D}$, i.e., when no enough data have been used in the screening procedure. Nonetheless, in this experiment, we have not encountered infeasible or suboptimal solutions.}

\subsubsection{Experiment 3 -- Worst Case}

 \textcolor{black}{Now} we evaluate the performance of the proposed method in an out-of-range situation. To do so, the time periods of the training and test sets are not randomly selected. Instead, we deliberately include the 1440 time periods with the highest aggregate demand values in the test set, and the remaining time periods in the training set. In this manner, we  artificially construct an unlikely situation that is, however, very adverse to the methods and, in particular, to UB+CC (which is the one that relies the most on the training data of the four constraint-screening approaches we consider in this experiment).

Similarly to Table \ref{tab:Large}, we provide the results of this out-of-range analysis in Table \ref{tab:Adverse}. As in the base case, the proposed UB+CC is the method that screens out the largest number of line-flow constraints by far, thus producing reduced UC problems that can be solved in around one fifth of the time that is needed to solve the full UC formulation. However, unlike in the base case, UB+CC leads to UC solutions that are infeasible or slightly suboptimal in a few instances. Interestingly, both UB and CC are able to deliver feasible and optimal UC plans in all instances of the test set under this very adverse situation. This reinforces the idea that the proposed cost-based valid inequality is only really effective, i.e., it has screening power,  when combined with a demand set $\mathcal{D}$ that does not contain values of the nodal net demands that are too different from the observed ones. This is so because, for a demand value that is too dissimilar from the observations, the upper bound may be too loose or even wrong. Therefore, the design of the set $\mathcal{D}$ is key to the performance of the proposed UB+CC, as this ``worst-case setup'' reveals.

\begin{table}[]
\caption{Results -- Base Case}
\centering
\begin{tabular}{lccccccc}
\hline
 & BN & UB & CC & UB+CC \\
\hline
Retained constraints (\%)   &33.8 &33.0 &24.9 & 18.9 \\
\#Infeasibilities   &0 &0 &0 &0 \\
\#Sub-optimal solutions   &0 &0 &0 &0 \\
Cost error (\%) &0.000 &0.000 &0.000 &0.000 \\
Reduced UC time (s)   &8.62 &8.15 &5.28 &3.42 \\
Computational burden (\%)   &74.6 &70.5 &45.7 &29.6\\
\hline
\end{tabular}
\label{tab:Large}
\end{table}

\begin{table}[]
\caption{Results -- Worst Case}
\centering
\begin{tabular}{lcccc}
\hline
 & BN & UB & CC & UB+CC \\
\hline
Retained constraints (\%)   &34.0 &31.9 &24.1 & 12.9 \\
\#Infeasibilities  &0 &0 &0 &22 \\
\#Sub-optimal solutions   &0 &0 &0 &55 \\
Cost error (\%)  &0.000 &0.000 &0.00 &0.008 \\
Reduced UC time (s)   &6.48 &5.99 &4.76 &1.86 \\
Computational burden (\%)   &61.2 &57.3 &45.5 &17.8\\
\hline
\end{tabular}
\label{tab:Adverse}
\end{table}

{\color{black} 
\subsubsection{\textcolor{black}{Experiment 4 -- Topological changes}}

We conclude this case study by assessing the resilience of the proposed constraint-screening approach, i.e., UB-CC, against topological changes. To this end, we have picked 34 corridor lines connecting automatic generation control areas of the 2000-bus system~\cite{2000bus} and assumed that \emph{any} of these lines (which are specified in~\cite{Screenpowersystems2021}), but only one at the same time, can be out of service. This leads to $34+1$ possible system configurations, where the ``$+1$'' corresponds to the base case analyzed and discussed in Experiment 2. 

Now, for each of those 35 possible topological configurations, we have run UB-CC following steps 1) and 2) of the procedure described in Section II-D, with each of these runs delivering a set of superfluous line-flow constraints. We have then computed the intersection set, that is, the set of the line-flow constraints that, according to UB-CC, can be safely removed from the original UC formulation under \emph{any} of those 35 possible configurations (these are the constraints that are detected as dispensable by UB-CC in all the 35 runs of the screening method).

As indicated in Table V, UB-CC identifies that 81.1\% of the line-flow constraints are superfluous in the base case with no topological changes. This percentage slightly drops (as should be expected) to 79.5\%, which shows that the screening rate of UB-CC remains high even under topological changes. This result suggests that our approach is resilient against this type of changes.
}

\section{Conclusions and Future Research}
\label{sec:conclusion}

The computational time of the unit commitment problem can be significantly reduced by screening out line capacity constraints. However, the optimization-based methods proposed in the literature only remove redundant constraints and thus, involve moderate computational savings. This paper presents a novel constraint screening methodology that removes both redundant and inactive constraints and further reduces the computational burden of this problem.

As existing approaches, the one we propose is based on computing the maximum line power flows on an LP-relaxation of the UC formulation. As a salient feature of our work, we propose to tighten this LP-relaxation to exclude uneconomical operating conditions. In doing so, our methodology is able to filter out a higher number of line capacity constraints. Simulation results using a 2000-bus network show that our proposal reduces the number of retained constraints and the solution time by 15\% and 45\%, respectively, if compared with existing benchmark methods. {\color{black} Furthermore, the constraint-screening rate of our approach remains quite unaltered when topological changes of the network are considered, which suggests that our approach is resilient against this type of changes. Finally, even though the cost inequality we use to increase the constraint-screening power of our method is data-driven, our numerical analysis reveals that the solution to the reduced UC problem we produce is generally feasible and optimal in the original UC formulation, if enough data are available.}


An aspect that requires further investigation is how the proposed methodology can be applied to multi-period unit commitment formulations that include inter-temporal constraints such as ramping limits and minimum times. In this case, the challenges are twofold, namely, i) building a sufficiently tight relaxation of the multi-period UC to produce an algorithm with a reasonable screening power, and ii) designing a proper set of net demands $\mathcal{D}$ and finding a tight upper bound on the system operating cost for a given time horizon.

%

\bibliographystyle{ieeetr}
\bibliography{mybibfile}
\begin{IEEEbiography}[{\includegraphics[width=1in,height=1.25in,clip,keepaspectratio]{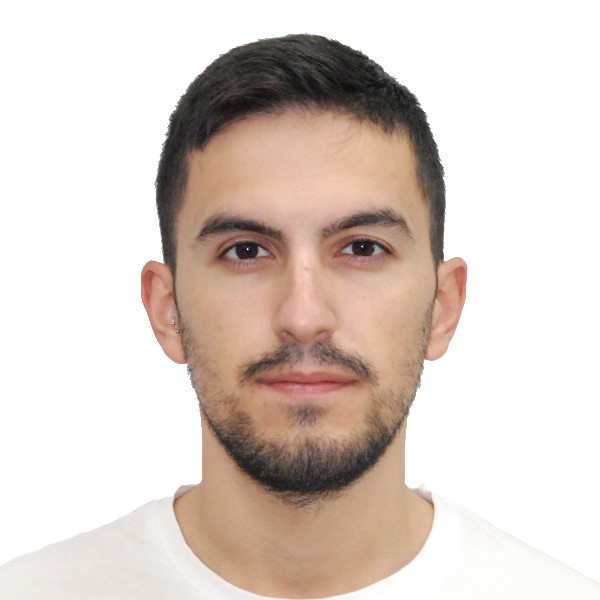}}]{\'Alvaro Porras} received the B.S. in Industrial Technologies Engineering and M.S. in Industrial Engineering from the University of M\'alaga, M\'alaga, Spain, in 2018 and 2020. He is currently a Ph.D. student at the University of M\'alaga and his research is financially supported by the Spanish Ministry of Science, Innovation and Universities through the university teacher training program.
\newline
\indent His research interests include power systems operation, decision-making under uncertainty and optimization.
\end{IEEEbiography}

\begin{IEEEbiography}[{\includegraphics[width=1in,height=1.25in,clip,keepaspectratio]{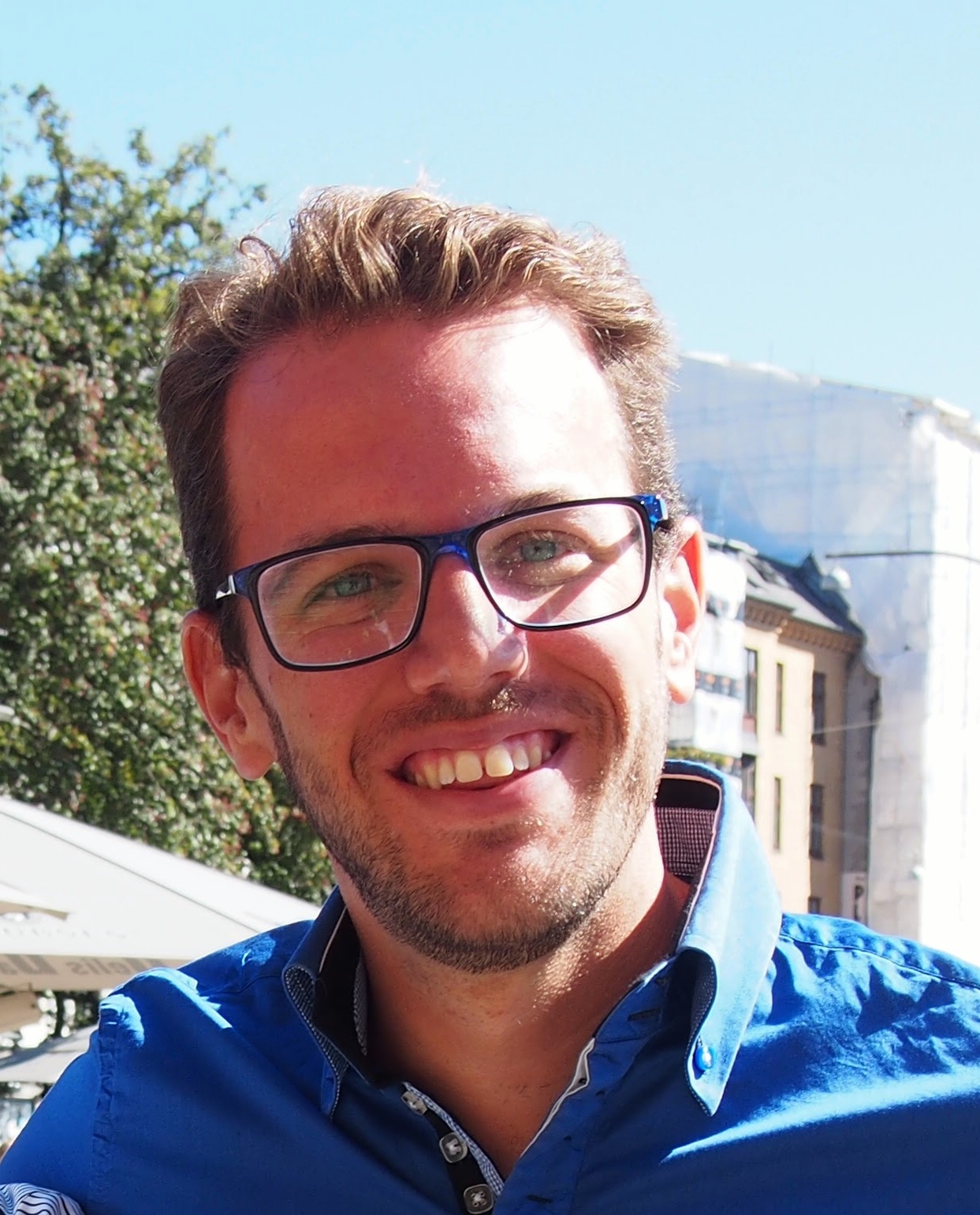}}]{Salvador Pineda} (S'07-M'11) received the Ingeniero Industrial degree from the University of  M\'alaga, Spain, in 2006, and a Ph.D. degree in Electrical Engineering from the University of Castilla-La Mancha, Spain, in 2011. He is currently an associate professor in the Department of Electrical Engineering at the University of M\'alaga in Spain.
\newline
\indent His research interests are in the fields of power system operation and planning, decision-making under uncertainty, bilevel programming, machine learning and statistics. 
\end{IEEEbiography}

\begin{IEEEbiography}[{\includegraphics[width=1in,height=1.25in,clip,keepaspectratio]{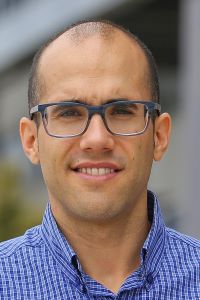}}]{Juan M.\ Morales} (S'07-M'11-SM'16) received the Ingeniero Industrial
degree from the University of  M\'alaga, M\'alaga, Spain, in 2006,
and a Ph.D. degree in Electrical Engineering from the University of Castilla-La Mancha, Ciudad Real, Spain, in 2010. He is
currently an associate professor in the Department of Applied
Mathematics at the University of M\'alaga in Spain.
\newline
\indent His research interests are in the fields of power systems
economics, operations and planning; energy analytics and optimization; smart
grids; decision-making under uncertainty, and electricity
markets.
\end{IEEEbiography}

\begin{IEEEbiography}[{\includegraphics[width=1in,height=1.25in,clip,keepaspectratio]{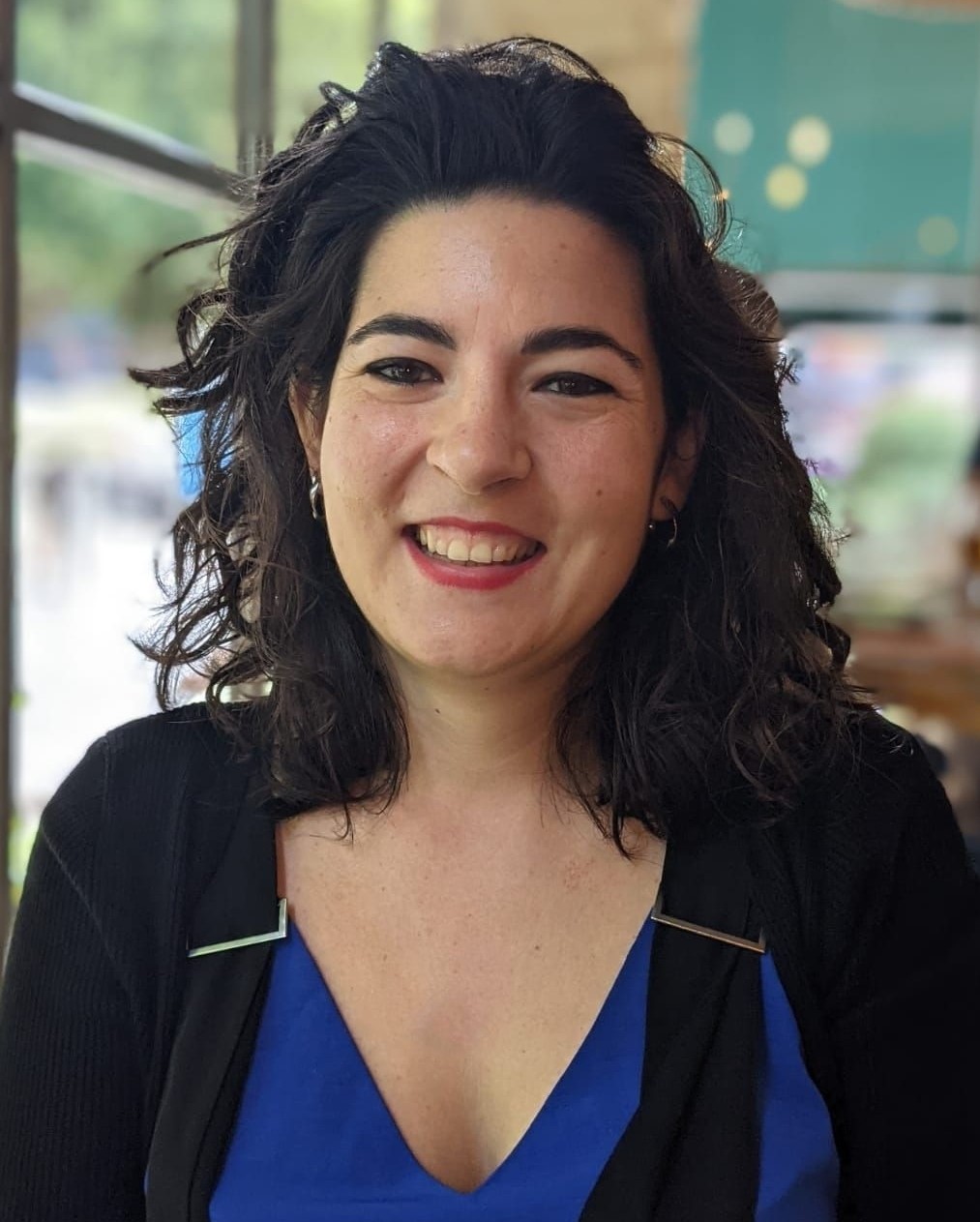}}]{Asunci\'on Jim\'enez-Cordero} received the Mathematics degree from the University of Seville, Seville, Spain, in 2013,
and a Ph.D. degree in Mathematics also from the University of Seville, in 2019. She is currently a lecturer in the Department of Statistics and Operations Research at the University of Málaga in Spain.
\newline
\indent Her research interests are in the fields of mathematical programming; optimization; machine learning; data-driven approaches, and power systems applications.
\end{IEEEbiography}

\end{document}